\input amstex
\documentstyle{amsppt}
\NoBlackBoxes
\vsize=8.9truein
\hsize=6.0truein

\def\rann{\operatorname{r\text{-}ann}}
\def\ZZ{{\Bbb Z}}
\def\deg{\operatorname{deg}}
\def\spec{\operatorname{spec}}
\def\grspec{\operatorname{gr\text{-}spec}}
\def\E{{\Cal E}}
\def\gr{{\operatorname{gr}}}

\def\Go{{\bf 1}}
\def\GL{{\bf 2}}
\def\LLV{{\bf 3}}
\def\NNV{{\bf 4}}
\def\NV{{\bf 5}}

\topmatter
 \title  The Graded Version of Goldie's Theorem \endtitle
\author K. R. Goodearl and J. T. Stafford\endauthor

\address Department of Mathematics, University of California, Santa
Barbara, Ca 93106-3080\endaddress
\email goodearl\@math.ucsb.edu \endemail
\address Department of Mathematics, University of Michigan, Ann Arbor,
MI 48109-1103 \endaddress
 \email jts\@math.lsa.umich.edu \endemail

\thanks The research of both authors was partially supported by grants
from the National Science Foundation. \endthanks



\abstract The analogue of Goldie's Theorem for prime rings is proved for
rings graded by abelian groups, eliminating unnecessary additional
hypotheses used in earlier versions. \endabstract

\endtopmatter  

\document

\head Introduction\endhead

In recent years, rings with a group-graded structure have become
increasingly important and, consequently, the graded analogues of
Goldie's Theorems have been widely utilized. Unfortunately, the graded
result requires an awkward extra condition: given a semiprime Goldie,
$\Bbb Z$-graded ring $R$, one cannot assert that $R$ has a
graded-semisimple ring of quotients unless one makes some extra
assumption, typically about the existence of homogeneous regular
elements (see \cite{\NV, Theorem~C.I.1.6}, for example), or about the
nondegeneracy of products of homogeneous elements (see
\cite{\NNV, Proposition~1.4}, for instance). The standard counter-example
\cite{\NV, Example~C.I.1.1} is the ring
$R=k[x]\oplus k[y]$, graded by giving $x$ degree $1$ and $y$ degree
$-1$. Note that this ring has no homogeneous regular elements other than
units,
 yet it is neither graded-semisimple nor graded-artinian.

There seems to be a misconception that Goldie's Theorem for prime rings
also requires such an extra condition (see, for example, \cite{\NV, 
\S{}C.I.1} or \cite{\LLV, p.~42}) 
and this is awkward in applications, as is illustrated by \cite{\GL,
 \S6.1} and \cite{\Go, \S5.4}. The purpose of this note is to correct this misconception by
showing that, at least for prime rings graded by abelian groups, 
no such extra hypotheses are required.

\head A Graded Goldie Theorem \endhead

Fix  an abelian group $G$.
Throughout,   ``graded'' will  mean ``$G$-graded''
 and, following standard practice,
  the graded analogue of a standard definition
 will be denoted by the prefix ``gr-''. Thus, for example, a gr-uniform 
 module  means a graded module that does not contain the direct sum of
 two nonzero, graded submodules. The formal definitions   can be found
in \cite{\NV}. 

The main aim of this note is to prove:

\proclaim{Theorem 1} Let $G$ be an abelian group and $R$ a 
$G$-graded, gr-prime, right gr-Goldie ring. Then, $R$ has a
gr-simple, gr-artinian right ring of fractions.
\endproclaim

In the conclusion of Theorem~1, it is tacitly assumed that the right
ring of fractions is taken with denominators which are homogeneous
regular elements of $R$, so that the $G$-grading on $R$ extends
(uniquely) to a $G$-grading on the ring of fractions.

The  only place where our proof differs significantly from
 its ungraded predecessors is in the last paragraph of the proof
of Theorem~4, where we need an extra trick to ensure that our chosen
regular element is also homogeneous. In the next three results, we only
require that $G$ be an abelian semigroup.

\proclaim{Lemma 2}  Keep the hypotheses of Theorem~1. Then:
\itemitem{\rm (i)} Any nonzero, graded right ideal $I$ of $R$ contains a
non-nilpotent homogeneous element. 
\itemitem{\rm (ii)}  The right gr-singular ideal of $R$ is nilpotent, and
hence zero. 
\endproclaim

\demo{Proof} Use the proof of
\cite{\NV, Lemma~C.I.1.4}, respectively \cite{\NV, Lemma~C.I.1.2}. (Those
results are stated for $\ZZ$-graded rings, but the structure of the
group does not enter into the proof.)
\qed\enddemo

\proclaim{Lemma 3}  Keep the hypotheses of Theorem~1. Suppose that
$a\in R$ is a homogeneous element  such that $aR$ is gr-uniform. Then,
its right annihilator,
$\rann(a)$, is maximal among right annihilators of nonzero homogeneous
elements of $R$. \endproclaim

\demo{Proof} 
 Suppose that $\rann(a) \subsetneq J= \rann(b)$, for some homogeneous 
element $b\in R$. Then $aJ\not=0$, and  so the gr-uniformity of
$aR$ implies that $aJ$ is gr-essential in $aR$. Therefore $aR/aJ$ is a
gr-singular right $R$-module. But $aR/aJ \cong R/J\cong bR$ because
$J\supseteq
\rann(a)$, and so $bR$ is gr-singular. Thus, by Lemma~1(ii), $b=0$.
\qed\enddemo

\proclaim{Theorem 4} Keep the hypotheses of Theorem~1.
 Then, any essential, graded right ideal $I$ of $R$ contains a
homogeneous regular element. \endproclaim

\demo{Proof} 
  Define a homogeneous element $a\in R$ to be  {\it gr-uniform} if  the 
right ideal $aR$ is gr-uniform. By Lemma~2(i) and the   Goldie
hypothesis, there exists a  non-nilpotent, gr-uniform element $a_1\in
I$. By induction, suppose that we have found  non-nilpotent, gr-uniform
elements $a_1,\dots,a_m\in I$ such that  $a_i\in 
\bigcap_{1\leq j\leq i-1} \rann(a_j)$ for $1<i\leq m$.  If 
$X=\bigcap_{1\leq j\leq m} \rann(a_j)\not=0$, then $I\cap X\not=0$ and 
so Lemma~2(i) again produces   a non-nilpotent, gr-uniform element
$a_{m+1}\in I\cap X$.  Since $a_i\in \rann(a_j)=\rann(a_j^2)$ for $i>j$, 
it is easy
to see that the  sum $\sum_{i\geq 1} a_iR$ is an internal direct sum.
Thus, as $R$ has finite gr-uniform dimension,  the inductive procedure
is finite. In other words, there exists 
 some index $n$ with $\bigcap_{i\le n} \rann(a_i)= 0$.

Since $R$ is gr-prime and the $a_i$ are not nilpotent,
 $a_1^2 Ra_2^2 R\cdots a_n^2R \ne 0$.  Thus, we may find homogeneous
elements $s_1,\dots, s_{n-1}\in R$ such that 
$a_1^2 s_1 a_2^2 s_2\cdots s_{n-1} a_n^2\not=0$. Then, by  Lemma~2(i),
there exists a homogeneous element $s_n$ such that 
$$c= a_1^2 s_1 a_2^2 s_2\cdots s_{n-1} a_n^2s_n$$ is not nilpotent. Set
$$d_i= \left(a_is_i a_{i+1}^2 s_{i+1} \cdots s_{n-1} a_n^2s_n\right)
\left( a_1^2s_1 \cdots s_{i-2} a_{i-1}^2 s_{i-1} a_i\right)$$ for
$i=1,\dots,n$. Note that the $d_i$ are subwords of $c^2$ and so  they
are non-zero. Therefore, Lemma~3 implies that
$\rann(d_i)= \rann(a_i)$, for each $i$. Moreover, the sum $\sum_{i=1}^n
d_iR$  is direct, because $d_iR \subseteq a_iR$. Hence,
$$\rann(d_1 +\dots+ d_n)= \bigcap_{i=1}^n \rann(d_i)  = \bigcap_{i=1}^n
\rann(a_i)=0.$$ Note that each $d_i$ is a reordering of the letters of
$c$ and so, as $G$ is abelian, $\deg(d_i)=
\deg(c)$. Therefore $d_1 +\dots+ d_n$ is a homogeneous regular element
lying in
$I$. \qed\enddemo

\demo{Proof of Theorem~1} The proof of Theorem~1 from Theorem~4  is
essentially the same as that in the ungraded case and is left to the
reader (see \cite{\NV, Proof of Theorem~C.I.1.6}).
\qed\enddemo

\head An Application\endhead

Recent motivation to prove Theorem~1 has come from the study of
quantized coordinate rings, where group-gradings can be used to
partition prime spectra in useful ways (see \cite{\GL, Section~6} and
\cite{\Go, Part~II} for a full discussion). The cited work required an
extra hypothesis, which can now be removed.

First consider an arbitrary group $G$ and a $G$-graded ring $R$, and let
$\grspec R$ denote the set of gr-prime ideals of $R$. For any ideal
$P$ of $R$, let $(P:\gr)$ denote the largest graded ideal contained in
$P$. Note that if $P$ is prime, then $(P:\gr)$ is gr-prime. Hence,
there is a partition
$$\spec R= \bigsqcup_{J\in \grspec R} \spec_J R,$$
where $\spec_J R= \{ P\in \spec R \mid (P:\gr)=J \}$. When $R$ is
noetherian and $G$ is free abelian of finite rank, each $\spec_J R$ is
homeomorphic to the scheme of irreducible subvarieties of an affine
algebraic variety, as follows.

\proclaim{Theorem 5} Let $R$ be a right noetherian ring graded by an
abelian group $G$, and let $J$ be a gr-prime ideal of $R$.
\itemitem{\rm (a)} The set $\E_J$ of homogeneous regular elements of
$R/J$ is a right denominator set, and the localization $R_J=
(R/J)[\E_J^{-1}]$ is a gr-simple, gr-artinian ring.
\itemitem{\rm (b)} The localization map
$R\rightarrow R/J\rightarrow R_J$
induces a Zariski-homeomorphism of $\spec_J R$ onto $\spec R_J$.
\itemitem{\rm (c)} Contraction and extension induce mutually inverse
Zariski-homeomorphisms between $\spec R_J$ and $\spec Z(R_J)$.
\itemitem{\rm (d)} If $G$ is free abelian of rank $r<\infty$, then
$Z(R_J)$ is a commutative Laurent polynomial ring over the field
$Z(R_J)_1$ in
$r$ or fewer indeterminates.
\endproclaim

\demo{Proof} Without loss of generality, $J=0$. Part (a) follows from
Theorem~1. In view of Theorem~4, all nonzero graded ideals of $R$ meet
$\E_0$, and so $\spec_0 R$ consists of the prime ideals disjoint from
$\E_0$. Hence, part (b) follows from standard localization theory. Parts
(c) and (d) follow from \cite{\Go, Corollary~4.3 and Lemma~4.1(d)}.
\qed\enddemo

In the applications to quantized coordinate rings, the grading arises
from a torus action, as follows. Let $R$ be a right noetherian algebra
over a field $k$, let $H= (k^\times)^r$ be an algebraic torus over $k$,
and suppose that we have a rational action of $H$ on $R$ by $k$-algebra
automorphisms. This implies that $R$ is graded by the character group
$G= \widehat H$, where the homogeneous component $R_g$ corresponding to
a character $g$ is just the $g$-eigenspace for the $H$-action on $R$.
(See \cite{\GL, \S6.1} or \cite{\Go, \S5.1} for details.) In this
setting,
$\grspec R$ is just the set of $H$-prime ideals of $R$, and $(P:\gr)=
\bigcap_{h\in H} h(P)$ for any ideal $P$. Theorem~5 applies, and the set
$\E_J$ is just the set of regular $H$-eigenvectors in $R/J$. The field
$Z(R_J)_1$ coincides with
$Z(R_J)^H$, the fixed subfield of $Z(R_J)$ under the induced $H$-action.
In fact, $Z(R_J)^H= Z(\operatorname{Fract} R/J)^H$, by the argument of
\cite{\Go, Proof of Theorem~5.3(c)}.

\enddocument